\documentclass{amsart}

\usepackage{amsmath}
\usepackage{amsthm}
\usepackage{amsfonts}
\usepackage{amssymb}
\usepackage{graphics}

\usepackage{xcolor}
\definecolor{refkey}{gray}{.45}
\definecolor{labelkey}{gray}{.45}

\newcommand\rank{\operatorname{rank}} 

\newcommand\R{{\mathbb{R}}}

\renewcommand\P{{\mathbf{P}}}
\newcommand\E{{\mathbf{E}}}

\usepackage{enumerate}
\usepackage{graphicx}
\usepackage[margin=1.35in]{geometry}
\usepackage[norelsize,ruled,vlined]{algorithm2e}

\newcommand{\N}{\mathbb{N}}

\newcommand\bu{{\mathbf u}}

\newcommand\BI{{\mathbf I}}

\newcommand\bv{\mathbf v}

\newcommand\CE{{\mathcal E}}
\newcommand\CF{{\mathcal F}}

\newcommand\CN{{\mathcal N}}

\newcommand\CS{{\mathcal S}}

\subjclass{26C10, 30C15}

\theoremstyle{plain}
  \newtheorem{theorem}{Theorem}

  \newtheorem{lemma}[theorem]{Lemma}
  \newtheorem{corollary}[theorem]{Corollary}
  \newtheorem{claim}[theorem]{Claim}

\theoremstyle{definition}

  \newtheorem{remark}[theorem]{Remark}

\begin{document}

\title{A simple SVD  algorithm for finding hidden partitions}

\author{Van Vu}
\address{Department of Mathematics, Yale, New Haven, CT 06520}
\email{van.vu@yale.edu}
\thanks{V. Vu is supported by research grants from NSF and Airforce.}

\begin{abstract} 
Finding a hidden partition in a random environment is a general and important problem, which contains as subproblems many famous questions, 
such as finding a hidden clique, finding a hidden coloring,  finding a hidden bipartition etc.

In this paper, we provide a simple  SVD algorithm for this purpose, answering a question of McSherry. This algorithm is very easy to implement and  works for sparse graphs with optimal density. 
 \end{abstract}

\maketitle

\section{The problem and a new algorithm } \label{section:introduction}

The hidden partition problem is the following: let $X$ be a set of $n$ vertices with a partition  $X = \cup _{i=1}^k X_i$;  for all $1 \le i \le j \le n$ and any $x \in X_i, y \in X_j$, we put a random edge between $x$ and $y$ with probability $p_{ij}$. Given one such random graph, the goal is   to recover the sets $X_i$.  This problem is of importance in computer science and statistics  and contains 
as special cases several well-studied problems such as  the  hidden clique,  hidden bisection, hidden coloring, clustering etc (see, for instance, 
 \cite{AK, AKS, AzarMc, BlumS, Bop, BCLS, CK, DGP, DM, FK, FR, JS, DF, Kucera, KVbook} and the references therein).  In what follows, we refer to  $X_i$ as clusters. 

\vskip2mm

In an  influential paper   \cite{MS},  Mc Sherry provided a  (randomized)  polynomial time algorithm that solves the general hidden partition 
problem for a large range of parameters. As corollary, he derived  several  earlier results  obtained  for   special cases.

\vskip2mm The general  idea \cite{MS} (and in many earlier works on clustering) is to find a good geometric representation of the vertices. We say that a representation is {\it perfect} if 
there is a number $r >0$ such that 

\begin{itemize} 
\item Vertices in the same cluster have distance at most $r$ from each other. 
\item Vertices from different clusters have distance at least $4r$ from each other. 
\end{itemize}

Once a perfect representation is obtained, it is easy to find the clusters. If $r$ is known, then the solution is obvious. If $r$ is not known, then there are several simple algorithms. For instance, one can create a minimal spanning tree (with respect to the distances) on the vertices and then remove the largest $k-1$ edges. In what follows, we put all these simple algorithms under a  subroutine
called  {\it Clustering by Distances} and the reader can choose his/her favorite to implement.  Our main goal is to present a simple way to obtain a perfect representation.

In the rest of the paper, let $s_u := |X_i|$ if $u \in X_i$ and $s := \min _{u \in X } s_u = \min_i |X_i| $.  We assume that $n$ is sufficiently large, whenever needed. Asymptotic notation 
are used under the assumption  $n \rightarrow \infty$. All explicit constants (such as  the $4$ above) are adhoc and we make no attempt to optimize them.

A popular way to find a perfect representation is to project the points  of $X$ (seen as vectors in $\R^n$) onto a properly chosen low-dimensional subspace $H$.
The  main technical part of Mc Sherry's  algorithm is a subroutine called $CProj$ (Combinatorial Projection), which 
creates  $H$ in a combinatorial way. The inputs  in this subroutine are a matrix  $\hat A$, parameters $k,s$, and a properly chosen threshold $\tau$. 


\begin{algorithm} \label{CPROJ} 
 (1) While there are at least $s/2$ unclassified nodes, choose an unclassified node $v_i$ randomly and define $T_i :=\{ u|  \| P_{\hat A^T } (\hat A^T_{v_i} -\hat A^T_u) \| \le \tau \}$. 
Mark each $u \in T_i$ as classified. 

(2) Assign each remaining node to the $T_i$ with the closest projected $v_i$.

(3) Let $\hat c_i$ be the characteristic vector of $T_i$. 

(4) Return $P_{\hat c} $, the projection onto the span of the $\hat c_i$. 
\caption{Combinatorial Projection (CProj)}  \end{algorithm}

\begin{algorithm}\label{algo:MS}
 (1)  Randomly partition the  set $\{1, \dots, n \}$ into two parts $A$ and $B$. Let $\hat A, \hat B$ be the submatrices of the adjacency matrix formed by columns from $A$ and $B$. 
 
 (2)  Let $P_1 = CProj (\hat B), P_2 =CProj (\hat A)$ and compute  $\bar H = [P_1 (\hat A) | P_2 (\hat B)]$. 
 
 (3) Run {\it Clustering by Distances} on the projected points. 
 
\caption{Mc Sherry's algorithm} 
\end{algorithm}

\vskip2mm


Let $P$ be the probability matrix  $(p_{ij})_{1 \le i, j \le n}$. For a vertex $u \in X$, $\bu$ denotes the corresponding column  in $P$. Define 

$$\Delta := \min  \| \bu  - \bv \| , $$ where the minimum is taken over all pairs $u, v$ belonging to different clusters. Mc Sherry proved \cite{MS}

 \begin{theorem}  \label{theorem:MS}  Assume that $\sigma^2  \gg \log^6 n /n $  is an upper bound on the variances of the entries. 
 There is a constant $C>0$ such that if 
 \begin{equation} \label{eqn:MScond}  \Delta   \ge  C    \sigma k^{1/2}   (\sqrt {\frac{n}{ s}}  +  \sqrt { \log \frac{n}{\epsilon} } ),  \end{equation} 
 the above algorithm  (with a proper choice of the threshold $\tau$)  recovers the partition with probability  $1-\epsilon$ with respect to the random graph and $k^{-1} $ with respect to the auxiliary random bits. 
  \end{theorem}

The main open question raised by Mc Sherry in \cite{MS} is  
to find  a more  natural and simpler  algorithm, which does not involve the subroutine  CPROJ  (see \cite[Section 4.4]{MS}).   The goal of this paper is to answer this question.

 To this end, $M_k$ denotes the subspace spanned by the first $k$ left singular vectors of 
 a matrix $M$.  Let $\hat P$ be our input, namely  the adjacency matrix of a random graph generated by $P$. 
 Arguably,   the most natural choice for $H$ would be   $\hat P_k$ (SVD), which leads to the algorithm below

\begin{algorithm}\label{algo:Vu1}
(1)    Project the columns of $\hat P$ onto $\hat P_k$.

 (2) Run {\it Clustering by distances } on the projected points. 
\caption{SVD  I}
\end{algorithm}

While SVD  I could well win the contest for being the simplest algorithm,  it is not easy 
to analyze in the general case.  In what follows, we   analyze a   slightly more technical  alternative,  SVD  II, which is a variant  of 
an algorithm proposed in \cite[Section 1]{MS}.


\begin{algorithm}\label{algo:Vu1}
 (0)  Randomly partition $X$ into two subsets $Y$ and $Z$. Let $B$ be the adjacency matrix of the bipartite graph between $Y$ and $Z$.   Let $Y_1$ be a  random subset of $Y$ by selecting each element with probability $1/2$ independently and  let $\hat A$ be the submatrix of $B$ formed by the columns indexed by $Y_1$. 
 
 (1) Project the columns of $B$ indexed by $Y_2:= Y \backslash Y_1$ on $\hat A_k$. 
 
 (2) Run {\it Clustering  by Distances } on the projected points. 
\caption{SVD  II}
\end{algorithm}

Compared to SVD  I, the extra steps in SVD  II  are the random partitions in Step $(0)$ done in order to reduce the correlation. 
   (A careful reading of \cite{MS} reveals that one also  need an extra partition in Algorithm \ref{algo:MS} to make the analysis go through.)

  Notice that SVD  II gives a partition  of $Y_2$, not $X$. There are many ways to extend it to a partition of $X$. For instance, we can  run the algorithm $l $ times (for some small $l$) and find partitions of 
$Y_2^{1}, \dots, Y_2^{l}$, where  $Y_2^i$ are random subsets of $X$ with density $1/4$ (the input graph is the same, only the random partitions are different). If a cluster $C$ in $Y_2^{i}$ and a cluster $C'$ in $Y_2^{i'}$ intersect, then they must belong to the same cluster in $X$ and we can merge them. If we choose $l=3 \log n$, say, then with probability $1-o(n^{-1} )$, all vertices of $X$ must belong to some $Y_2^{i}$ and we recover the clusters $X_1, \dots, X_k$ at the end. 
We can also first  find the partitions of $Y_1, Y_2$ and $Z_1, Z_2$ by reversing the role of $Y_1$ and $Y_2$ and $Y$ and $Z$ and find which four clusters must belong to an original cluster by looking at the edge densities; we omit the details. 

Beside being simple, SVD  II is also very convenient to implement, as its main step, the computation of the  projection onto $\hat A_k$ (given $\hat A$ as input) is a routine operation (SVD) which appears in most standard 
mathematical packages.

 Let us now analyze SVD  II.   For convenience, we assume that $P$ has rank $k$. 
 The general case when $P$ can have a smaller rank is  discussed later.   Let  $\lambda$ be the least non-trivial singular value of $P$.

 \begin{theorem} \label{theorem:main1}  There is a constant $C>0$ such that the following holds. Assume that   $\sigma^2  \ge C \frac{\log n}{n}$ and $s \ge C \log n,  k = o( (n/\log n)^{1/2} )$. 
Then  SVD  II
  clusters $Y_2$ correctly  with probability  $1 - o(n^{-1}) $ if one of the following two conditions is satisfied

\begin{itemize}  
\item   Condition 1.    $\Delta \ge C ( \sigma \sqrt {\frac{n}{s}} + \sqrt {\log n } ).$ 
 
 \item  Condition 2. 
     $\Delta  \ge C ( \sigma \sqrt {\frac{n}{s} }  +   \sigma \sqrt{ k \log n} +  \frac{ \sigma \sqrt {nk}} { \lambda}  )$ 
 \end{itemize}

 If we omit the assumption $s \ge C \log n$, the statement still holds but with probability $$1 -o(n^{-1} ) -  c\sum_{i=1}^k e^{-|X_i| /c}  $$ for some constant $c$.
  \end{theorem}

\begin{remark} We would like to point out a few remarks 
\vskip2mm 

\begin{itemize}

\item The lower bound  $\sigma^2 \ge C \log n/n$ is optimal, up to the value of $C$. If $\sigma^2 <  \log n/n$, then there are many isolated points, which can be assigned to any cluster.

\item We can reduce the  failure probability $O(n^{-1} )$ to  $O(n^{-K }) $ for any constant $K$ at the cost of  increasing the constant $C$.

\end{itemize} 
\end{remark}

\vskip2mm Let us now consider the performance of SVD  II on various subproblems.  We allow the value of $C$ to be flexible in order to omit smaller order terms  for convenience. 
It is instructive  to compare the corollaries below with Corollaries 1,2,3 from \cite{MS}. 

{\it Hidden clique.} In this problem,  $k=2$, $s$ is the size of the clique, and $\Delta = (1-p) \sqrt s$, where $p$ is the density of the random graph. 
 Condition 1 becomes 

$$(1-p) s^{1/2}  \ge C ( p^{1/2}  \sqrt {\frac{n}{s}  }  + \sqrt {\log n}  ) $$ which is satisfied if $s \ge C ( \sqrt {np}   + \sqrt {\log n}  ) $.  As $np  = \Theta ( \sigma^2 n ) = \Theta (\log  n)$, this simplifies to 
$s \ge C \sqrt {np} $.  

\begin{corollary} There is a constant $C$ such that for any $p \ge C \frac{\log n}{n}$ and $s \ge C \sqrt {np}$,  SVD  II finds the hidden clique of size $s$ with probability $1- o(1)$. 
\end{corollary}

{\it Hidden Coloring.} Here $k$ is the number of color classes,  each has size $n/k$;  $\Delta = p \sqrt {2n/k} $; $s =n/k$; $\sigma^2 = p(1-p) $.    The singular values of $P$ are $\frac{k-1} {k} n , \frac{1}{k} n, \dots, \frac{1}{k} n $. 
If $p \ge 1/k$, Condition 1 is
$$p \sqrt {n/k}  \ge C ( p^{1/2} \sqrt { k} + \sqrt {  \log n } )  $$ which is satisfied for  $k = o( (n /\log n))^{1/3} )$. 

 If $ p < 1/k$, then the bound  $\lambda \ge \sigma \sqrt {ns}$ holds, and the $\Delta$ bound in Condition 2 is $$p \sqrt {n/k} \ge C (\sqrt {p k \log n} + \frac{k \log n}{\sqrt n } )$$  which is satisfied if $p \ge C \frac{ k^{3/2} \log n }{n} $.

\begin{corollary} There is a constant $C$ such that the following holds.
For any $k =o( (n/\log n)^{1/3}$ and  edge density $.99 > p   \ge  C   \frac{ k^{3/2}  \log n}{n} $,   SVD  II finds the hidden $k$-coloring  with probability $1 - O(n^{-1} )$. 
\end{corollary} 

{\it  Hidden Bipartition.} Let the two densities be $.99 \ge p > q>0 $. We have $k=2$,  $\Delta = |p-q| n^{1/2} $, $s =n/2$, $\sigma^2 = \Theta (p) $.   The two singular values of 
$P$ are $(p+q)n$ and $(p-q)n$.  Condition 2 requires 
$\frac{p-q}{p^{1/4} }  \ge C  \sqrt { \frac{\log n }{  n } }. $

\begin{corollary} There is a constant $C$ such that the following holds
Let $.99 > p >q   \ge C \log n /n $ be edge densities such that  $\frac{p-q}{p^{1/4} }  \ge C \sqrt {\frac{\log n}{n } } $ then SVD  II finds the hidden bipartition with probability $1 - o(n^{-1} )$. 
\end{corollary} 

One can replace $p^{1/4} $ in the denominator by a better term $p^{1/2} $ by considering  an approximate  algorithm; see  Corollary \ref{cor4}.


The rest of the paper is organized as follows. In the next section, we present a few technical lemmas and  prove Theorem \ref{theorem:main1}  and Theorem \ref{theorem:approx} 
in Section \ref{section:proof}.  In Section \ref{section:variants}, we discuss  variants of SVD  II,  including an approximate version
which works under weaker assumptions.

\section{Technical  lemmas}  \label{section:lemmas}

\begin{lemma} [Projection of a Random Vector]  \label{lemma:projection} There are constants $C_1, C_2$ such that the following holds. 
Let $\xi = (\xi_1, \dots, \xi_n)$ be a random vector in $\R^n$  whose coordinates $\xi_i$ are independent  random variables with 
 mean 0 and variance at most $\sigma^2 \le 1$. Let $H$ be a subspace of dimension $d$ and $\Pi_H \xi$ be the length of the orthogonal projection 
of $\xi$ onto $H$. Then

$$\P ( \Pi_H X  \ge   \sigma \sqrt{d} +   C_1 \sqrt {\log n} ) \le n^{-3} . $$ 

Furthermore, if $H$ has an orthornormal bases  $v_1, \dots, v_d$ such that  $\max_{ 1\le i \le d } \| v_i \| _{\infty} \le \alpha$, then 

$$\P ( \Pi_H X \ge C_2 \sqrt{d} (\sigma + \alpha \log n ) )  \le n^{-3} . $$ 
\end{lemma}

We prove this lemma in the appendix.



\begin{lemma} [Norm of a random matrix] \label{lemma:norm} There is a constant $C_0 >0$ such that the following holds. 
Let $E$ be a symmetric matrix whose upper diagonal entries $e_{ij}$  are independent random variables where $e_{ij} = 1- p_{ij}$ or 
$-p_{ij}$ with probabilities $p_{ij}$ and $1-p_{ij}$, respectively, where $0 \le p_{ij} \le 1$. Let $\sigma^2 := \max _{ij} p_{ij}(1-p_{ij}$.  
 If  $\sigma^2 \ge C_0  \log n / n $, then 

$$ \P ( \| E   \|  \ge  C_0  \sigma n^{1/2} ) \le  n^{-3}. $$

\end{lemma} 

If $\sigma^2 \ge \frac {\log^4 n}{n} $, the statement  is a corollary of  \cite[Theorem 1.4]{Vunorm}. 
For smaller $\sigma$, one can  prove this lemma using the $\epsilon$-net approach by Kahn and Szemeredi \cite{KS}. We omit the details, which is very similar to the proof of 
Feige and Ofek  for \cite[Theorem 1.1]{FO}.


\begin{lemma}  [Perturbation bound]  \label{lemma:DK}  Let $M, N$ be  matrices where $\delta := \lambda _k (M) - \lambda_{k+1} M >0$. Then 
$$\sin \angle (M_k, (M+N)_k )  \le \delta^{-1}  \| N \| . $$
\end{lemma}

This lemma is  a well known result in numerical linear algebra, known as Davis-Kahan-Wedin theorem; see \cite{B, DK, W, GVL}.

\section{Proof of Theorems \ref{theorem:main1}  }  \label{section:proof}

  Let  $A$ be the probability matrix $p_{ij}$ corresponding to $\hat A$.  As $A$ is a large random submatrix of $P$, it is not hard to show  that $\lambda_k (A) =\Theta ( \lambda_k (P) ) $ with high probability (we provide a verification of this fact at the end of the proof). In the rest of this proof, we assume 

\begin{equation} \label{eqn:lambda}  \lambda_k (A) \ge c_0 \sigma \sqrt {ns} , \end{equation}  for some constant $c_0 >0$.

We  view the adjacency matrix  $\hat A$ (between $Y_1$ and $Z$)  as a random perturbation of $A$,  $\hat A := A+E $, where the  entries $e_{ij}$ of $E$ 
are independent and $e_{ij} = 1 - p_{ij}$ with probability $p_{ij}$ and $ -p_{ij}$ with probability $1-p_{ij}$.  We denote by $\hat \bu, \bu, e_u$ the columns corresponding to a vertex 
$u$ in $\hat A, A, E$, respectively.  All matrices are of size approximately $n/2 \times n/4$ by the definitions of $Y,Z$ and $Y_1, Y_2$. 

Our leading idea is that the random perturbation $E$ does not change $A_k$ too much, thus hopefully the projections onto $\hat A_k$ and $A_k$ differ 
by only a small amount. The heart of the matter, of course, is to bound this error term. While inviting, a straightforward  application of  Lemma \ref{lemma:DK} is 
too crude in the general case (it does lead to some simple solution for some subproblems in certain range of parameters). We will still make use of  this lemma, but for a
quite  different purpose.

For simplicity, we assume in the rest of the proof that $s \ge C \log n$. For a sufficiently large $C$, this  implies that with probability $1 -o(n^{-1} )$, each cluster $X_i$ intersects $Z$ in at least  $|X_i|/3$ elements. Thus, the distance between two columns (belonging to different clusters) in $A$ is at least $\Delta /3$.  We aim 
to show that with high probability $\| P_{\hat A_k } \hat \bu - \bu \| < \Delta/12$ for all $u \in Y_2$; this will provide  a perfect geometric representation.   If there is no lower bound on $s$, then the probability that the random partition has this property is at least $1 - c\sum_{i=1}^k e^{-|X_i| /c} $ for some constant $c>0$.

For a fixed $u$, by the triangle inequality

$$\| P_{\hat A_k } \hat \bu - \bu \|  \le \| P_{\hat A_k}  (\hat  \bu - \bu)\|  + \| (P_{\hat A_k}  - I ) \bu \| = \| P_{\hat A_k}  e_u \|  + \| (P_{\hat A_k}  - I ) \bu \| . $$

To bound the second term, we follow an argument from \cite{MS} and consider 

$$(P_{\hat A_k} -I ) A = (P_{\hat A_k } - I) \hat A -  (P_{\hat A_k } -I) E . $$

The spectral norm of the first term is $\lambda_{k+1} (A_k) \le \lambda_{k+1} (A ) + \| E\| = \|E \|$, as $A$ has rank at most $k$. 
The spectral norm of the second term is also at most $\| E \|$. Thus, by Lemma \ref{lemma:norm}, by probability at least 
$1 - n^{-3}  $ 

$$\| (P_{\hat A_k } - I ) A \| \le 2 \| E \| \le C_0 \sigma n^{1/2} , $$  for some constant $C_0$. 

Let $\chi_u$ be the unit vector $s_u^{-1/2} \BI_u$ where $\BI_u$ is the indicator vector for the cluster containing $u$,  we have 

$$\|( P_{\hat A_k } -I ) A \| \ge \| (P_{\hat A_k } -I ) A \chi_u  \|  = s_u^{1/2} \| (P_{\hat A_k } -I ) \bu \| . $$

Combining the last two inequalities and using the union bound, we conclude that with probability at least $1 - n^{-2} $ 

$$\| (P_{\hat A_k}  - I ) u \| \le C_0  \sigma \sqrt {\frac {n}{s_u} } ,  $$ for all $u \in X$.

Now we tend to the first term, whose analysis is more involved.   By the first part  of Lemma \ref{lemma:projection}, 

$$\| \ P_{\hat A_k } e_u \| \le \sigma k^{1/2} + C_1 \sqrt { \log n  } $$ with probability $1 - o(n^{-2} )$, for a properly chosen constant $C_1$.  As $sk \le n$, the  term $\sigma k^{1/2} $ is at most 
$\sigma \sqrt {n/s}$ and  can be omitted. This yields  that if 
 
$$\Delta \ge  C_0  \sigma \sqrt {n/s} +  C_1 \sqrt {\log n  }$$  then the algorithm  succeeds with probability at least $1 - o(n^{-1} ) $. 
This  proves the first part of the theorem concerning Condition 1.

To prove the second part of the theorem,  let us reconsider  the distance 
$P_{\hat A_k } e_u $.   Notice that if $s \le 10 k  \log n $, then Condition 2 implies Condition 1 (with some  modification on the value of $C$). Thus,  in what follows, we can assume $s \ge 10 k \log n $.

Rewrite $\hat A = A +  E$ and let $v$ be a singular vector of $A$. Recall that $|X_i \cap Z| \ge  \frac{1}{3} |X_i| = s_i/3$ for all $i$ 
By symmetry, each  coordinate in $v$ is repeated at least $s/3$ times, thus $\| v\| _{\infty} \le  2s^{-1/2} $. 
Furthermore, by  Lemma \ref{lemma:DK} and Lemma \ref{lemma:norm}, we have with probability $1 - o(n^{-2})$ that  

$$\sin (A_k, \hat   A_k ) \le  C_0 \frac{\sigma \sqrt n } {\lambda}  $$

\noindent which implies that  for any unit vector $v \in \hat A_k$, $$\| v\|_{\infty}  \le  2 s^{-1/2} +  C_0 \frac{\sigma \sqrt n } {\lambda}  \le C_0' s^{-1/2} $$ by the condition on 
$\lambda$, with some properly chosen constant $C_0'$.   Using the second part of Lemma \ref{lemma:projection}, we conclude that with probability $1 -o(n^{-2} )$,
$\| P_{\hat A_k } e_u  \|  \le C( \sigma k^{1/2} +  s^{-1/2} \log n )$ for all $u$ and some properly chosen constant $C$, concluding the proof.

For the sake of completeness, let us show  that with high probability, the least (non-trivial) singular value of  $P$ and $A$ are essentially the same, up to a constant factor. 
We first compare the singular values of $P$ with the singular values of  $\tilde P$, the probability matrix of the bipartite graph spanned by $Y$ and $X$. Using Chernoff's bound, one  can easily show that 
with probability at least $1 -n^{-2} $ 

\begin{equation}  \label{intersection}  | |X_i \cap Y | - |X_i|/2|   \le 5 \sqrt {|X_i| \log n } \end{equation}  for all $1 \le i \le k$. 

We use the fact that for any matrix $M$ of rank $k$
$\lambda _{k}  (M) = \inf_{rank (M') = k-1}  \| M - M' \| _F $. For simplicity, let us assume for a moment that $|X_i \cap Y| = |X_i|/2$. Let $\tilde P'$ be the matrix that define $\lambda_k (\tilde P) $. 
We define $P'$, a rank $(k-1)$ approximation of $P$, by extending $\tilde P'$ as follows. For the block indexed by $X_i \backslash Y$, simply copy the block of $\tilde P'$ corresponding to $X_i \cap Y$. It is trivial that $P'$ has  rank $k-1$ and 

$$\| P -P'\| _F^2 =  2 \| \tilde P -\tilde P' \| _F^2 $$ which implies $\lambda_k  \le \sqrt {2} \lambda_k (\tilde P ) $. With the same argument, we can compare $\lambda_k (\tilde P)$ with $\lambda_k (B)$ and the later with $\lambda_k (A)$, each time losing a factor of $\sqrt 2$. At the end it would give  $\lambda_k (P) \le 2^{3/2} \lambda_k (A)$. 

To make the argument precise, we need to remove the assumption $|X_i \cap Y| = |X_i| /2$. Using \eqref{intersection}, we can create a matrix $P'$ such that 

$$\| P-  P'\|_F^2 \le 2 \| \tilde P -\tilde P'\| _F ^2 + 5 \sum_{i=1}^k \sqrt{|X_i| \log n} \sigma ^4 . $$ 

On the other hand, the extra term $ 5 \sum_{i=1}^k \sqrt{|X_i| \log n} \sigma ^4 $ is less than $\frac{1}{4} \lambda_k (P)^2$ by the assumption of the theorem.  Thus, we can use the above estimate to get a slightly weaker bound 
$\lambda_k (P) \le 2  \lambda_k (\tilde P) $, completing the proof.

\section {Variants }  \label{section:variants}

\subsection{Dimension and Density}  
In the case $\rank A = l < k$, it makes more sense to project onto $\hat A_l$ rather than onto $\hat A_k$; the rest of the algorithm and analysis remains the same. 

If we do not  know either $k$ or $l$ in advance. 
we can   modify SVD  II slightly as follows.  The idea is to define the essential rank of $\hat A$ to be the largest index $l$ such that 
$\lambda_l (\hat A ) \ge C_3 \sigma n^{1/2}$, for a properly chosen $C_3$, 
 and replace the projection onto $\hat A_k$ by the projection onto $\hat A_l$.  After  redefining $\lambda := \lambda_l P $, the content  of Theorem \ref{theorem:main1} remains the same. Its proof also remains  the same, expect few nominal changes. The error term caused by smaller singular values is  not going to effect the final conclusion.

The only information we need in SVD  II (using essential dimension) is the value of $\sigma$.  Even if this information is not known,
 we can still solve the problem by considering a sequence of $O(\log n)$  trials with $\sigma_1  = \frac{\log n }{\sqrt n }, \sigma_i =2 \sigma_{i-1} $ and run SVD  II in each case. Each trial will output a clustering and it is easy to decide which  one is correct   by considering the degree densities of the vertices from  one cluster  to the others.

\subsection{ An approximate Solution.} 
In practice, one is often  satisfied with an approximate solution.
We say that   a partition $X =\cup_{i=1}^k  X_i'$ is $\epsilon$-correct if $|X_i \backslash X_i' | \le \epsilon |X_i| $.  
Similarly, we say that a geometric representation of $X$  is $\epsilon$-perfect  if 
there are points $x_1, \dots, x_k$ with distance at least $4r$ from each other so that at least $(1-\epsilon)|X_i|$ points from $X_i$ has distance at most $r$ to $x_i$. 
One can use an $\epsilon$-perfect representation to find an $\epsilon$-correct partition.

 
 \begin{theorem} \label{theorem:approx}  Given $\epsilon >0$, there  is a constant $C>0$ such that the following holds.   If  $\sigma^2  \ge C  \frac{\log n}{n}$ 
 and  $$\Delta  \ge C \sigma \sqrt {\frac{n}{s} } ,$$ then 
with probability $1-o(n^{-1} )$  the projection in SVD  II produce an $(1-\epsilon)$-perfect representation of the point sin $Y_2$.  \end{theorem}

It is worth  mentioning   that in various situations, an $\epsilon$-correct partition can be upgraded to a fully correct one by a simple ``correction" procedure, as shown in the following example: 

{\it Hidden bipartition.} Assume $p > q$.  
Let $X= X_1' \cup X_2 '$ be an $\epsilon$-correct partition, for some small $\epsilon$ (say $\epsilon=.1$).  Then both $X_i'$ have size at most $\frac{1}{2} (1-\epsilon ) n $.  Assume $|X_i \backslash X_i'| \le \epsilon n/2$; it follows that 
$|X_1' \backslash X_1 | \le \epsilon n$. 
  With probability $1 - n^{-2} $, the following holds. For any $u \in X$, let $d_u$ be the number of its neighbors in $X_1'$. 
  If $u \in X_1$, then 
  
  $$ d_u  \ge |X_1 \cap X_1'| p + |X_1' \backslash X_1 | q - 10 \sqrt {np \log n }= D_1 . $$
  
  On the other hand, if $u \in X_2$, then 
  
  $$d_u \le |X_1 \cap X_1 '| q + |X_1' \backslash X_1 | p + 10 \sqrt {np \log n} = D_2 . $$

  It is clear that if $(p-q) \ge 30 \sqrt{ \frac {p \log n} {n} }$, then $D_1 > D_2$. Thus, one can correct the partition by defining $X_1$ be the set of $n/2$ vertices $u$ with largest $d_u$.

\begin{corollary} \label{cor4}  There is a constant $C$ such that the following holds
Let $.99 > p >q   \ge C \log n /n $ be edge densities such that  $\frac{p-q}{p^{1/2} }  \ge C \sqrt {\frac{\log n}{n } } $ then the approximation algorithm with correction
 finds the hidden bipartition with probability $1 - o(n^{-1} )$. 
\end{corollary} 


{\it Proof of Theorem \ref{theorem:approx}. } 
We first bound  $\| P_{\hat A_k } e_u \|$. Recall that 

$$\E \| P_{\hat A_k } e_u \| ^2 \le  \sigma^2 k .$$

By Markov's  inequality, it follows that $ \P ( \| P_{\hat A_k } e_u \|  \ge K \sigma k^{1/2} ) \le  K^{-2} $.  We call a vertex $u$ {\it good} if   $\| P_{\hat A_k } e_u \|  \le  K \sigma k^{1/2} $. 
For a sufficiently large $C$ (depending on $K$),  all good vertices will be clustered correctly. Moreover, choosing  $K \ge 2 \epsilon ^{-1/2}$, the probability for  $u$ being good 
is at least $1-\epsilon/4$, thus the expectation of the number of good elements in $X_i$ is at least $|X_i| (1 -\epsilon /4)$. As the good events are independent, Chernoff's bound implies that  with probability $1- n^{-2} $,  at least $|X_i| (1- \epsilon)$ points from $X_i$ are good. This completes the proof.

\appendix 

\section{Proof of Lemma \ref{lemma:projection}} 

Notice that the function 
$\Pi_H (X)$ is $1$-Lipschitz and convex, thus by Talagrand's inequality \cite{Tconc}   for any $t >0$

$$\P ( \Pi _H X \ge \mu + t) \le 2 \exp (-t^2/4) $$ where $\mu$  is the mean of $\Pi_H (X)$. We do not know $\mu$; however, we can bound from above. 
Slightly abusing the notation, let $\Pi := (\pi_{ij})  $ denote the projection matrix onto $H$, then 

$$\E  |\Pi_H X |^2  = \E X^T \Pi X = \sum_{i=1}^n \pi_{ii} \E \xi_i^2 \le  \sigma ^2 \sum_{i=1}^n \pi_{ii} = d \sigma^2. $$

Combining this with the concentration inequality, it is not hard to show that $\mu \le \sigma d^{1/2} + O(1) $, concluding the proof of the first part of the lemma.

The second part  follows immediately  from 

\begin{claim} \label{Chernoff} Let $(a_1, \dots, a_n)$ be real numbers such that $\sum_i a_i^2 =1$ and $|a_i | \le \alpha$ for all $i$. Let 
$\xi_i$ be independent random variables with mean 0 and $\E |\xi_i| ^k \le \sigma^2$ for all $k \ge 2$. Let $S:=\sum_{i=1}^n a_i \xi_i$. Then 

$$\P ( |S| \ge   4 (\sigma \sqrt {\log n } +    \alpha \log n )   \le 2n^{-3} . $$  \end{claim}

To prove Claim \ref{Chernoff}, notice that   for any $0 < t \le \alpha^{-1} $ we  have 

$$\E \exp (t S )  = \prod_i  \E \exp( t a_i \xi_i ) = \prod _i (1+  \frac{ \sigma^2 a_i^2  t^2 }{2!} + \frac{ t^3 a_i^3 \E \xi_i^3 }{6 ! } + \dots  )  $$

Since $\E \xi_i^k \le \sigma^2$ for all $k \ge 2$ and $t|a_i | \le 1$, the right most formula is

$$ \le \prod_i (1+ \sigma^2 t^2 a_i^2) \le \exp (\sigma^2 t^2 ). $$

Markov's inequality yields 

$$\P ( S \ge  T ) \le \exp( -t T + t^2 \sigma^2 ) . $$

To optimize the RHS, let us consider two cases 

{\it Case 1.}  $\sigma \ge \alpha  \sqrt { \log n}  $. Take $T = 4 \sigma  \sqrt {\log n} $ and $t =  \frac{ \sqrt {\log n}} { \sigma }  \le \alpha^{-1} $. 
With this setting  $-tT + t^2 \sigma^2 = - 3  \log n$.

{\it Case 2.} $\sigma <  \alpha \sqrt { \log n}  $. Take $T = 4 \alpha \log n $  and $t = \alpha^{-1} $.  In this setting,   $-tT + t^2 \sigma^2  \le -4 \log n + \log n =- 3 \log n $. 

One can bound $\P ( -S \le T)$ the same way.

\vskip2mm 

{\it Acknowledgement.} The author would like to thank NSF and AFORS for their support and K. Luh for his careful proof reading. 





\end{document}